\documentclass[12pt,leqno]{article}
\usepackage{amsmath,amssymb,amsfonts}
\usepackage{eucal}

\numberwithin{equation}{section}

\newtheorem{lem}{Lemma}[section]

\newtheorem{prop}[lem]{Proposition}
\newtheorem{thm}[lem]{Theorem}
\newenvironment{proof}[1][Proof]{\textbf{#1.} }{\ \rule{0.5em}{0.5em} \medskip}

\newtheorem{exam}{Example}[section]

\renewcommand{\epsilon}{\varepsilon}
\renewcommand\tilde{\widetilde}

\newcommand\inv{^{-1}}
\renewcommand\Im{\operatorname{Im}}     

\newcommand\Ker{\operatorname{Ker}}

\newcommand\FM{\operatorname{FM}}
\newcommand\G{\Gamma}
\newcommand\GG{(\G,g,G)}

\newcommand\cF{\mathcal{F}}
\newcommand\cG{\mathcal{G}}
\newcommand\cT{\mathcal{T}}

\newcommand\bbZ{\mathbb{Z}}

\begin{document}

\title{Cycle and Circle Tests of Balance in Gain Graphs: \\
    Forbidden Minors and Their Groups 
    }
\author{Konstantin Rybnikov\footnote{Department of Mathematical Sciences, University of Massachusetts at Lowell, Lowell, MA 01854, U.S.A.  Part of this research was conducted at Cornell University, Ithaca, NY 14853-4201, U.S.A.}  \\
    and Thomas Zaslavsky\footnote{Department of Mathematical Sciences, Binghamton University, Binghamton, NY 13902-6000, U.S.A.  Research partially supported by National Science Foundation grant DMS-0070729.} 
    }

\maketitle

\abstract{We examine two criteria for balance of a gain graph, one based on binary cycles and one on circles.  
The graphs for which each criterion is valid depend on the set of allowed gain groups.  The binary cycle test is invalid, except for forests, if any possible gain group has an element of odd order.  Assuming all groups are allowed, or all abelian groups, or merely the cyclic group of order 3, we characterize, both constructively and by forbidden minors, the graphs for which the circle test is valid.  It turns out that these three classes of groups have the same set of forbidden minors.  The exact reason for the importance of the ternary cyclic group is not clear.}

\bigskip



\bigskip\hrule\bigskip

\emph{Note to publisher:}

This paper does NOT have a ``corresponding author''.

All authors are EQUAL.

All authors are able to answer correspondence from readers.

\bigskip

For editorial purposes ONLY, contact the writer of the cover letter of submission.

\bigskip\hrule\bigskip


\section{Introduction}\label{intro}

A \emph{gain graph} $\GG$ is a graph $\G=(V,E)$ together with a group $G$, the \emph{gain group}, and a homomorphism $g$, the \emph{gain function}, from the free group $F(E)$ on the edge set $E$ into $G$.  We think of the edges of $G$ as oriented in an arbitrary but fixed way, so that if $e$ is an edge in one direction, then $e\inv$ is the same edge in the opposite direction; thus $g(e\inv) = g(e)\inv$.  A gain graph is \emph{balanced} if every simple closed walk lies in the kernel of $g$.  One of the fundamental questions about a gain graph is whether or not it is balanced.  We examine two related approaches to this question.

A simple, general test for balance is to examine the gains of a fundamental system of circles.
A \emph{circle} is the edge set of a nontrivial simple closed walk (a walk in which no vertex or edge is repeated, except that the initial and final vertices are the same).  If we take a spanning tree $T$ of $\G$ (without loss of generality $\GG$ can be assumed connected), each edge $e \notin T$ belongs to a unique circle in $T\cup e$.  These circles constitute the \emph{fundamental system of circles} with respect to $T$.  The simple closed walk corresponding to a circle $C$ is unique up to choice of initial vertex and direction.  Thus it depends only on $C$, not on the choice of walk, whether the walk is in $\Ker g$.  If it is, we say $C$ is \emph{balanced}.  
It is well known and easy to prove that $\GG$ is balanced if and only if every circle of a fundamental system (with respect to some spanning tree) is balanced (see, for instance, the generalization in
Zaslavsky (1989), Corollary 3.2).

The trouble with testing a fundamental system of circles is that what one knows about the gains may not be about a fundamental system.  
Thus we look for a more general sufficient condition that a gain graph be balanced.  
A fundamental system of circles is one kind of basis of the binary cycle space $Z_1(\G;\bbZ_2)$; 
we generalize by considering an arbitrary basis (giving what we call the Binary Cycle Test for balance) or an arbitrary basis composed of circles (giving the Circle Test).  
Our study is based on the fact that, given a basis that is not a fundamental system, one cannot always decide balance of the gain graph by testing the basis; it may be impossible to reach a decision.  
(We might remind the reader here that since a closed walk in a graph, taken as a binary cycle, is reduced modulo $2$, a walk $n$ times around a circle $C$ gives the binary cycle that is $C$ itself if $n$ is odd but the zero cycle if $n$ is even.  A walk that goes around a circle $C_1$ once, then takes a path to another circle $C_2$ disjoint from $C_1$, goes around $C_2$ once, and returns on the same path, reduces to the disconnected binary cycle $C_1 \cup C_2$: the path disappears when taken modulo $2$.)
  
In this paper we take an arbitrary class of gain groups and focus on a particular underlying graph $\G$,\footnote{In Rybnikov and Zaslavsky (20xx) we consider arbitrary graphs, but only with abelian gain groups.} asking whether, for every basis $B$ of the appropriate type (depending on the test being used) and every gain mapping $g$ into a particular gain group (or any of a specified list of gain groups), $\GG$ is necessarily balanced if $B$ lies in the kernel of $g$.  (This description is vague in several ways.  We shall make it precise in the next section.)  If this is so, we call $\G$ \emph{good}.  (Again, a precise definition will follow.)  We find four principal results.

First, the class of good graphs is closed under deletion and contraction of edges---that is, any \emph{minor} (a contraction of a subgraph) of a good graph is again good (Theorem \ref{T:minors}).  This is true for any choice of allowed gain groups and for all the tests that we define.  It follows that the good graphs are characterized by a list of forbidden minors, and by the main theorem of Robertson and Seymour's ``Graph minors'' series 
(1985, 20xx),
this list is finite.  The natural problem is to find the forbidden minors explicitly.  In a broader sense, one wants to know whether there are many graphs for which each test is valid. This, of course, depends on which test is being applied and what is the list of permitted gain groups. 

Our second main result defines the range of validity of the Binary Cycle Test by finding the forbidden minor (Theorem \ref{T:fmcycle}).  If an odd-order group is a possible gain group then the Binary Cycle Test works only for forests (Theorem \ref{T:fmcycle}); thus, it is useless, since a forest is always balanced.  However, if the gain group has no elements of odd order, the Binary Cycle Test may be of use; deciding for which abelian gain groups that is so is the topic of 
Rybnikov and Zaslavsky (20xx)
(see Theorem \ref{T:abelian} below), while for nonabelian groups this question is open.  
In another direction, we may take Theorem \ref{T:fmcycle} as a suggestion to restrict the construction used to calculate gains (see the end of Section \ref{cycle}).

Our third result is curious (Theorem \ref{T:good}).  For the Circle Test there are exactly four forbidden minors and they are bad graphs if and only if amongst the possible gain groups is the cyclic group of order 3.  
Thus $\bbZ_3$ seems to have special importance.  This fact is surprising and mysterious.  Might it be related to the fact that all forbidden minors are small?

The even wheels and even double circles suggest in a qualitative way that this may be true, since $W_{2k}$ and $2C_{2k}$ are bad if $\bbZ_{2k-1}$ is a gain group (our fourth significant result, Theorem \ref{L:fmw4}); hence, for any gain group with odd torsion there is a bad graph.  
Moreover, the larger the wheel or double circle, the more are the groups for which it is bad.  However, it lies beyond the power of our methods to explain the way in which the class of admissible gain groups influences the class of good graphs.



\section{Definitions}\label{prelim}

The graphs that we consider may be infinite and may have loops and multiple edges.  (A \emph{loop} in a graph is an edge whose two endpoints coincide.  A non-loop edge is called a \emph{link}.)  
A \emph{closed walk} is a sequence of vertices and edges, $v_0e_1v_1e_2 \cdots e_lv_l$, in which the endpoints of $e_i$ are $v_{i-1}$ and $v_i$, starting from and ending at the same vertex $v_0=v_l$.  It is \emph{trivial} if $l=0$.  
It is \emph{simple} if it is nontrivial and it does not repeat any vertex or edge except for having the same initial and final vertex.  
A \emph{theta graph} is a graph homeomorphic to a triple link.

(Gain graphs have been called \emph{voltage graphs} in the context of surface embedding theory, where the actual gain around a closed walk is important.  See, e.g., 
Gross and Tucker (1987).  
For us it only matters whether the gain is the identity.  We eschew the term ``voltage'' because gains do not have to obey Kirchhoff's voltage law.)

\emph{Switching} a gain graph $\GG$ means replacing $g$ by $g^f$, obtained in the following way: take any function $f: V \to G$ and for an edge $e$ oriented with initial vertex $v$ and final vertex $w$, define $g^f(e) = f(v)\inv g(e) f(w)$.  
Switching does not change which circles are balanced, nor whether the gain graph is balanced. 
It is easy to see that $G$ can be switched so that, in any chosen maximal forest $T$, every edge has identity gain: $g|_T \equiv$ the group identity.

A \emph{binary cycle} is the indicator function of a finite edge set that has even degree at every vertex; thus, we may identify the group of binary cycles with the class of all such edge sets, with symmetric difference as the addition operation. 
We write $Z_1(\G;\bbZ_2)$ for the group of binary cycles, where $\bbZ_2 = \bbZ/2\bbZ$.  
(In topological language, a binary cycle is a 1-cycle in the cellular homology of $\G$ with coefficients in $\bbZ_2$.)

Suppose that $b$ is a binary cycle.  A \emph{cyclic orientation} of $b$ is any closed walk $\tilde b$ whose abelianization taken modulo 2, that is, whose natural projection into $Z_1(\G;\bbZ_2)$, is $b$.  
If $B$ is a set of binary cycles, a \emph{cyclic orientation of $B$} is any set $\tilde B = \{ \tilde b : b \in B \}$ of cyclic orientations of the members of $B$.  
A cyclic orientation is necessarily connected, but a binary cycle with disconnected support can still have a cyclic orientation (if $\G$ is connected) because an edge in $\tilde b$ need not be in the support of $b$.  
There is no cyclic orientation of $b$ that can reasonably be regarded as canonical, except in the case of a binary cycle whose support is a circle.

We can now state the Binary Cycle Test.

\medskip
\emph{Definition.}   
Let $\GG$ be a gain graph and $B$ be a basis of $Z_1(\G;\bbZ_2)$.  We say that $B$ \emph{passes the Binary Cycle Test} if it has a cyclic orientation $\tilde B$ such that all elements of $\tilde B$ have gain 1.  
We say the Binary Cycle Test is \emph{valid for $\GG$} if the existence of a basis $B$ that passes the binary cycle test implies that $\GG$ is balanced.  
We say the Binary Cycle Test is \emph{valid} for a family of graphs $\cF$ and a family of groups $\cG$ if it is valid for every gain graph $\GG$ with $\G \in \cF$ and $G \in \cG$.
\medskip

In other words, the Binary Cycle Test is valid for a gain graph $\GG$ if the existence of a basis $B$ of $Z_1(\G;\bbZ_2)$ with a cyclic orientation $\tilde B$, all whose members have gain 1, implies that $\GG$ is balanced. The converse implication is always true: if $\GG$ is balanced then every cyclic orientation of every binary cycle has necessarily gain $1$. 

Two examples, that will be used again later in proofs, will show that the Binary Cycle Test can indeed fail.

\begin{exam}\label{X:loop} {\rm 
$K^\circ_1$ is the graph consisting of a loop and its supporting vertex.
Suppose $k$ is odd and $k\geq 3$.  Take $\G = K^\circ_1$ with binary cycle basis $\{e\}$, and let the cyclic orientation of $e\in Z_1(\G)$ be $\tilde e = ee\cdots e$ ($k$ times).  Take gain group $\bbZ_k$ and assign gain $g(e) =$ a generator of $\bbZ_k$.  Then $g(\tilde e) =$ the identity, so the basis passes the Binary Cycle Test, but $(\G,g,\bbZ_k)$ is unbalanced.  Therefore the Binary Cycle Test is invalid for $(\G,g,\bbZ_k)$.
}
\end{exam}

\begin{figure}
\vskip 2.6in{}
\includegraphics{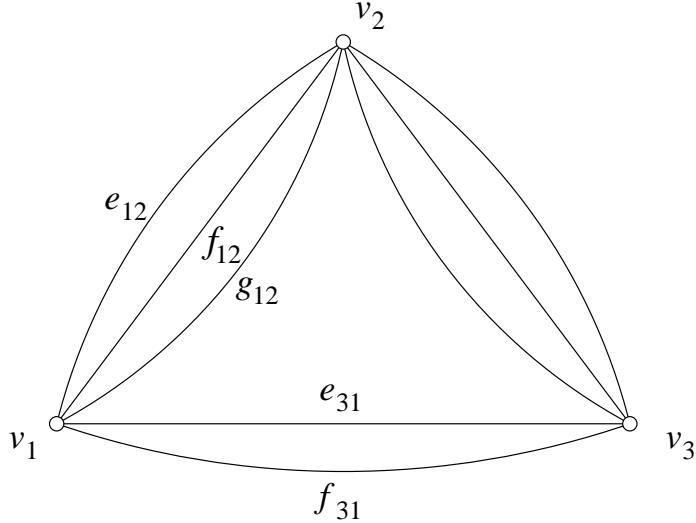}
\caption{$C_3(3,3,2)$ for Example \ref{F:c3(332)} and Proposition \ref{P:c3}.} \label{F:c3(332)}
\end{figure}

\begin{exam}\label{X:c3(332)} {\rm 
The graph $C_3(3,3,2)$ is shown in Figure \ref{F:c3(332)}.  
We show that $C_3(3,3,2)$ is bad if and only if $\bbZ_3 \in \cG$.  
A balanced circle basis $B$ that satisfies the digon property contains six triangles, no two having more than one common edge.  Without loss of generality we may take 
\begin{equation*}
B = \{ e_{12}e_{23}e_{31}, f_{12}g_{23}e_{31}, g_{12}f_{23}e_{31}, f_{12}f_{23}f_{31}, e_{12}g_{23}f_{31}, g_{12}e_{23}f_{31} \}
\end{equation*} 
and, by switching, we may take $g(e_{ij}) = 1$.  It follows that $g(f_{i-1,i}) = g(g_{i,i-1}) = a$, say, where $a^3 = 1$.  Thus if $\bbZ_3 \in \cG$, the Circle Test can fail; but otherwise it cannot. 
}
\end{exam}

We may wish to restrict the bases and cyclic orientations.
For instance, one often wants to apply the Binary Cycle Test only to cycles that are circles, with the natural cyclic orientation as a simple closed walk.  The \emph{Circle Test} is the Binary Cycle Test with those restrictions.  
Example \ref{X:c3(332)} shows that the Circle Test can fail.

Ideally one would wish to find all pairs $(\G,G)$ consisting of a graph and a group such that each of these tests is valid for all gain graphs $\GG$ (that is, for all gain mappings $g: E\to G$).  
We cannot give a complete solution to this problem, but, as outlined in the introduction, we do have a partial answer.


\section{Minor Closure}\label{minors}

A class of graphs (that is, of isomorphism types, or unlabelled graphs) is \emph{minor closed} if, for any graph in the class, all its minors are in the class.
A \emph{minor} of a graph $\G$ is the result of any finite sequence of successive operations of contracting edge sets and taking subgraphs.  If $\G$ is infinite, the edge set contracted and the subgraph taken may be finite or infinite.  It is easy to see that if the union of all contracted edge sets is $S$, then the minor is a subgraph of $\G/S$.  Thus only one operation of each type is required.  

\begin{thm}\label{T:minors}
Let $\cG$ be any class of groups.  The class of graphs $\G$ such that the Binary Cycle Test is valid for any gain graph $\GG$ with $\G$ as underlying graph and gain group in $\cG$ is minor closed.

The same holds for the Circle Test.
\end{thm}

We prove the theorem for a class $\cG$ consisting of a single group $G$.  Clearly, this implies the whole theorem.  
A subgraph can be taken by deleting an edge set and then removing any subset of the isolated vertices.  Since in connection with gains any isolated vertices are immaterial, it suffices to treat just deletion of edge sets and contraction.
The proof, therefore, consists of one lemma for deletion and one for contraction, and a remark on the Circle Test.

\begin{lem}\label{L:deletion}
If the Binary Cycle Test is valid for $\G$ and $S \subseteq E(\G)$, then it is valid for $\G \setminus S$.
\end{lem}

\begin{proof}
We have to prove that the Binary Cycle Test is valid in $\G \setminus S$.  That is, let $B$ be any binary cycle basis of $\G \setminus S$ and $\tilde B$ any cyclic orientation of $B$; let $g$ be any gain map on $\G \setminus S$ with gain group $G$; we must prove that $g$ is balanced.

We extend $g$ to gains on $\G$ in the following way.  First, let $T$ be a 
maximal forest of $\G/(E \setminus S)$ and let $S' = S \setminus T$.  
(That every infinite graph has a maximal forest, which is the union of a spanning tree in each component, was proved by K\"onig \cite[Chapter IV, Theorems 24 and 27]{K}.)  
Next, for each $e \in S'$, let $C_e$ be a circle of $(\G \setminus S') \cup e$ that contains $e$.  Let $\tilde C_e$ be a natural cyclic orientation.  Now we assign gain $1$ to each $e \in T$ and we assign gains to $e \in S'$ so that $g(C_e) = 1$.

Now, $B' = B \cup \{C_e : e \in S'\}$ is a basis for $Z_1(\G;\bbZ_2)$.  To see why, note that the edges of $T$ serve to connect components of $\G \setminus S$ that are connected in $\G$, while each edge $e$ of $S'$ increases the dimension of the binary cycle space.  
Since $Z_1(\G \setminus S,\bbZ_2)$ and the $C_e$ for $e\in S'$ span $Z_1(\G;\bbZ_2)$, if $B$ spans $Z_1(\G \setminus S,\bbZ_2)$ then $B'$ spans $Z_1(\G;\bbZ_2)$.
Since $C_e$ is the unique element of $B'$ that contains $e$, if $B$ is independent then so is $B'$.  
As we assumed $B$ to be a basis of $Z_1(\G \setminus S,\bbZ_2)$, $B'$ is a basis of $Z_1(\G;\bbZ_2)$.

Let $\tilde B'$ be a cyclic orientation of $B'$ extending that of $B$.  Since $B'$ is a cycle basis of $\G$ and $g(\tilde b) = 1$ for every $b \in B'$, the extended $g$ is balanced.
\end{proof}

\begin{lem}\label{L:contraction}
If the Binary Cycle Test is valid for $\G$ and $S \subseteq E(\G)$, then it is valid for $\G/S$.
\end{lem}

\begin{proof}
We have to prove the Binary Cycle Test is valid in $\G/S$.  Let $T$ be a maximal forest in $S$ and $S' = S\setminus T$.  Since $\G/S = (\G/T) \setminus S'$, by Lemma \ref{L:deletion} it suffices to show validity in $\G/T$.

Let $\cT$ be the set of connected components of $(V,T)$.  $\G/T$ has vertex set $\{ v_T : T \in \cT \}$ and edge set $E \setminus T$.  If the endpoints of $e \in E \setminus T$ are $u$ and $v$ in $\G$, its endpoints in $\G/T$ are $T_u$ and $T_v$, where $T_u$ denotes the component of $(V,T)$ that contains $u$.

We have a fixed binary cycle basis $B$ of $\G/T$ and cyclic orientation $\tilde B$.  We have to convert them to a binary cycle basis $B'$ of $\G$ and a cyclic orientation $\tilde B'$.  We first convert each $\tilde b \in\tilde B$.  Suppose $\tilde b = (T_0,e_1,T_1,\ldots,e_l,T_l)$ (the $T_i$ being vertices of $\G/T$ and $T_0 = T_l$), and each $e_i$ has endpoints $u_i \in T_{i-1}$ and $v_i \in T_i$ and is oriented so its direction in $\tilde b$ is from $u_i$ to $v_i$.  Let $T_{v_iu_{i+1}}$ be the unique path in $T_i$ from $v_i$ to $u_{i+1}$.  (We take subscripts modulo $l$.)  Then 
$$
\tilde b' = (v_0, T_{v_0u_1}, u_1, e_1, v_1, T_{v_1u_2}, u_2, \ldots, v_{l-1}, T_{v_{l-1}u_l}, u_l, e_l, v_l )
$$
and $b'$ is the projection of $\tilde b'$ into the binary cycle space $Z_1(\G;\bbZ_2)$.

 We should verify that $B'$ is a basis of $Z_1(\G;\bbZ_2)$; but this is obvious because $Z_1(\G;\bbZ_2)$ is naturally isomorphic to $Z_1(\G/T,\bbZ_2)$ and $B'$ naturally maps to $B$.

Now take $g$ to be any gain function of $\G/T$ with gain group $G$.  Then $g$ is also a gain function on $\G \setminus T$.  We extend $g$ to $\G$ by setting $g(e)=1$ for $e \in T$.  Thus, $g(\tilde b') = 1$ for every $b' \in B'$.  By the Binary Cycle Test in $\G$, $\GG$ is balanced; therefore $(\G/T,g,G)$ is balanced.
\end{proof}

\bigskip
\begin{proof}[Proof of Theorem \ref{T:minors} completed]  
The proof demonstrates minor closure of graphs that satisfy the Circle Test because the new cycles in Lemma \ref{L:deletion} are circles and the modified cycles in Lemma \ref{L:contraction} convert circles to circles.
\end{proof}

The theorem implies that for each test and every class $\cG$ of groups there is a list $\FM(\cG)$ of forbidden minors, finite graphs such that the test (whichever test it is) is invalid for some $G\in\cG$ and some $G$-gain graph based on each forbidden graph, but if a finite graph $\G'$ has none of the forbidden graphs as a minor, then for every $G\in \cG$ and every $G$-gain graph based on $\G'$ the test is valid.
Thus we have the following problem:  Given a class $\cG$ of groups, we consider all graphs $\G$ such that the Binary Cycle Test, or the Circle Test, is valid for all gain graphs $\GG$ with $G\in \cG$.  Of course, $\cG$ is a class of isomorphism types.  We should assume that it is \emph{subgroup closed}, that is, if $H < G\in \cG$, then $H\in \cG$, since we can produce all gain graphs with gains in $H$ simply by taking $\Im g$ in an appropriate subgroup of $G$.  We define $\FM_0(\cG)$ to be the class of forbidden minors for validity of the Circle Test for gain graphs with gain group in $\cG$ and $\FM_2(\cG)$ to be the corresponding class for the Binary Cycle Test.


\section{Homeomorphism, Whitney Operations, and Extrusion}\label{whitney}

Homeomorphic graphs are equivalent for our purposes.  (Two graphs are \emph{homeomorphic} if they are both obtained by subdividing edges of the same graph.  \emph{Subdividing an edge} means replacing it by a path of positive length; of course, if the length is 1 the subdivision is trivial.)

\begin{lem}  \label{L:homeo}  
Given any class $\cG$ and either of our two tests, if the test is valid for $\G$ then it is valid for any graph homeomorphic to $\G$.
\end{lem} 

\begin{proof}
We may suppose $\G$ has no divalent vertices.  Let $\G'$ be homeomorphic to $\G$ and let $g'$ be a gain function on $\G'$.  We may assume by switching that $g'$ is the identity on all but one edge $e_P$ of each maximal induced path $P$ of $\G'$.  If we contract all but that one edge in each $P$, we have $\G$ with a gain function that we call $g$.  It is clear that the correspondence $P \mapsto e_P$ defines bijections of the binary cycle spaces and of the cyclic orientations of binary cycles, as well as a gain-preserving bijection of closed walks, between $\G'$ and $\G$.  Thus binary cycle and circle bases of the two graphs correspond, and the kernels of $g'$ and $g$ correspond.  It is now easy to see that each test is valid in $\G'$ if and only if it is valid in $\G$.
\end{proof}

\emph{Whitney operations} on a graph are:  identifying two vertices in different components, the inverse of that operation, and twisting one half of a 2-separation.  Only the latter concerns us, since balance is a property that depends only on the blocks of $\G$.  A more precise definition of a twist in a block graph $\G$ is this: find a separating vertex set $\{u,v\}$ and for $i = 1,2,\ldots,k$ let $A_i$ be a component of $\G\setminus \{u,v\}$ together with $u$, $v$, and the edges connecting them to the component.  Let $A$ be any union of some but not all subgraphs $A_i$.  \emph{Twisting $A$ across $\{u,v\}$} means reconnecting to $v$ every edge of $A$ that was incident to $u$, and vice versa.

\begin{lem}  \label{T:twist}  
If the Circle Test is valid for a graph $\G$, it is valid for any graph obtained by twisting $\G$ across a 2-separation.   
\end{lem}

\begin{proof}  
Whitney operations do not alter the binary cycle space.  We borrow from 
Zaslavsky (2003), Section 5, 
the observation that, suitably interpreted, they preserve gains.  
The trick is that, when $A$ is twisted, the gain on $e\in A$ is reversed, so that $g'(e) = g(e)^{-1}$.  It is now clear that, if a closed walk had identity gain before twisting, it remains so after twisting.  
Thus a basis $B$ that had a cyclic orientation with identity gain before twisting continues to have such an orientation after twisting.  Also, the gain graph is balanced after the twist if and only if it was already balanced.
\end{proof}

A difficulty with Whitney twisting in connection with the Binary Cycle Test is that a particular cyclic orientation of a binary cycle may, after twisting, become disconnected.

\begin{figure}
\vskip 1.25in{}
\includegraphics{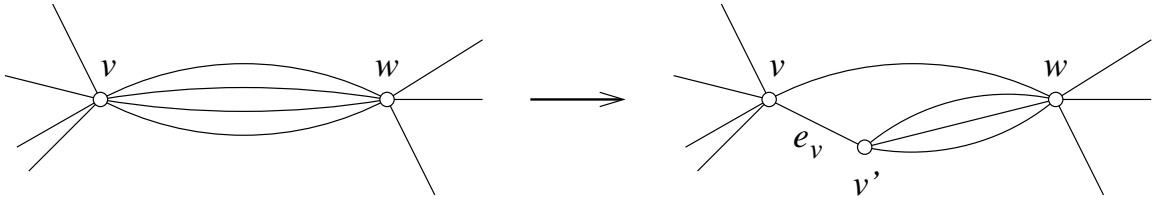}
\caption{Extruding an edge from $v$.} \label{F:extrude}
\end{figure}

The inverse operations to edge deletion and contraction do not in general preserve goodness.  The general inverse of contraction is called \emph{splitting a vertex}: $v$ is replaced by $v_1$ and $v_2$ and a new edge $v_1v_2$, and each edge incident with $v$ becomes incident with $v_1$ or $v_2$.  If we restrict splitting so that $v_2$ has degree two, the result is homeomorphic to the original graph so goodness is preserved.  This can be generalized.  \emph{Extruding} a vertex from $v$ means choosing a neighbor $w$ and at least one of the $vw$ edges and adding a new vertex $v'$ and edge $e_v$ so that $e_v$ joins $v$ to $v'$ and all the selected $vw$ edges become $v'w$ edges.  No other edges are affected.  (See Figure \ref{F:extrude}.)

\begin{lem} \label{L:extrude}
Given a subgraph-closed class $\cG$ of groups and a loopless graph $\G$ for which the Circle Test is valid, extruding an edge in $\G$ maintains the validity of the test.
\end{lem}

\begin{proof}
Suppose $v'$ is extruded from $v$ as in Figure \ref{F:extrude}, forming $\G'$.  We show that the circles in $\G'$ are essentially the same as those in $\G$.  Let $e_1,\ldots,e_k$ be the $v'w$ edges in $\G'$.  If a circle $C$ in $\G$ passes through just one of the edges $e_i$, then in $\G'$, $C$ will have a break $vv'$ that can be filled in by $e_v$.  Thus $C$ in $\G$, regarded as an edge set, becomes $C' = C\cup e_v$ in $\G'$.  Any other circle of $\G$ is a circle in $\G'$.  Conversely, a circle $C'$ in $\G'$ that does not contain $e_v$ cannot visit both $v$ and $v'$, because if it visits $v'$ it is a digon with vertex set $\{v',w\}$.  Thus, $C'$ is a circle in $\G$.

This analysis shows that circle bases of $\G$ and $\G'$ correspond and so do the exact expressions for circles in terms of a given basis.  Gains also correspond up to switching since $e_v$ can always be given gain $1$ while the gains of all other edges are the same in $\G$ and $\G'$.
\end{proof}

This lemma implies some but not all cases of the preceding one.  It applies only to the Circle Test and to single-edge extrusions, so it does not apply to a graph that is homeomorphic to $\G$ by the introduction of infinitely many divalent vertices.

A graph is \emph{extrusion-irreducible} if it is not obtainable by extrusion from a smaller graph.  An equivalent property, if the graph is loopless, is that a vertex with exactly two neighbors is multiply adjacent to both.

\begin{lem} \label{L:fmextrude}
Let $M_1,M_2,\ldots$ be a list of finite graphs with the properties that 
\begin{enumerate}
\item[(i)]  every vertex in $M_i$ has at least two neighbors, 
\item[(ii)]  $M_i$ has no loop, and 
\item[(iii)]  $M_i$ is extrusion-irreducible.
\end{enumerate}  If $\Gamma$ is a finite graph that has no minor isomorphic to any $M_i$, then extruding a vertex in $\Gamma$ yields a graph with no minor isomorphic to any $M_i$.
\end{lem}

\begin{proof}  
Suppose extruding $v'$ and $e_v$ from $v$ in $\Gamma$ makes a graph $\Gamma'$ with a subgraph $\Gamma_0$ that contracts to $M_i$.  Take $\Gamma_0$ to be minimal.  Then $v'$ must be in $\Gamma_0$, or else $\Gamma_0 \subseteq \Gamma$.  Also, $e_v$ must be an edge in $\Gamma_0$, or else $v'$ would have at most one neighbor in $\Gamma_0$ and, by (i), $\Gamma_0$ could not be minimal.  Further,
$e_v$ cannot be contracted in forming $M_i$ from $\Gamma_0$ because $M_i$ is not a minor of $\Gamma$.
At least one of the $v'w$ edges $e_1,e_2,\ldots$ of $\Gamma'$ must be in $\Gamma_0$ to avoid violating (i).  But then (iii) is violated unless one of them is contracted in forming $M_i$.  Say $e_1,\ldots,e_k$ are in $\Gamma_0$ and $e_1$ is contracted.  If $k > 1$ we get a loop in $M_i$ (because $\Gamma_0$ is minimal), contradicting (ii).  If $k = 1$, then contracting $e_1$ is equivalent to contracting $e_v$, which is impossible.  Therefore no $\Gamma_0$ can exist.  
\end{proof}

Another general result on extrusion is Lemma \ref{L:type1}.


\section{Forbidden Minor for the Binary Cycle Test}\label{cycle}

The result for the Binary Cycle Test is very simple.  
Recall that $K^\circ_1$ is the graph consisting of a loop and its supporting vertex.

\begin{thm}  \label{T:fmcycle}  
Let $\cG$ be the class of all groups, or any class containing a
nontrivial group of odd order.  Then $\FM_2(\cG) = \{K^\circ_1\}$.  
\end{thm}  

\begin{proof}  Suppose $\bbZ_k \in \cG$ where $k$ is odd and $k\geq 3$.  See Example \ref{X:loop} for a proof that there is a gain graph based on $K^\circ_1$ for which the Binary Cycle Test is invalid.
\end{proof}

What happens when $\cG$ contains only groups without odd-order elements is not known in general.  However, if the groups are all abelian we have results from 
Rybnikov and Zaslavsky (20xx):

\begin{thm} \label{T:abelian}
The binary cycle (or circle) test is valid for $\G$ if $\cG$ is a class of abelian groups without odd torsion and either $\G$ is finite or no group in $\cG$ has an infinitely 2-divisible element other than zero.
\end{thm}

For example, $\cG$ may consist of the trivial group and all groups having exponent 2, as all such groups are abelian.  
Gain groups of the form $\bbZ^r_2$ (whose gain graphs may be called \emph{multisigned graphs}) seem to be important, e.g., for nonorientable surface embedding.  Signed graphs themselves (group $\bbZ_2$) were treated in 
Zaslavsky (1981), 
whose argunents, although stated for circles, apply equally to binary cycles.  
Clearly, groups with exponent 2 are a very special case, for then there is nothing gained by going outside the binary cycle space to either cyclic orientations or integral cycles.


\section{Validity of the Circle Test}\label{circle}

In this section's main result we characterize the graphs for which the Circle Test is valid, so long as some allowed gain group has an element of order 3.

Throughout this section we consider only finite graphs.  $\cG$ denotes a class of groups that is closed under taking subgroups, and all gain graphs have gain group in $\cG$.

\emph{Definition.}
When considering a specific class of groups $\cG$, we call a graph $\G$ \emph{good} if the Circle Test is valid for all gain graphs $\GG$ with $G\in \cG$.

Some important graphs: 
\begin{itemize}
\item $W_n$ is the $n$-spoke wheel; its rim vertices in cyclic order are $v_1, v_2, \ldots, v_n$, the hub is $w$, and the rim is the circle $v_1 v_2 \cdots v_n v_1$.  
\item $m\G$ is $\G$ with every edge replaced by $m$ copies of itself. 
\item $C_l(m_1,\ldots,m_l)$ is a circle $C_l = e_1\cdots e_l$ where $e_i$ is replaced by $m_i$ copies of itself; thus $C_l(m,\hdots,m) = mC_l$.  
\item $K_4(m,m')$ is $K_4$ with two opposite edges replaced by $m$ and $m'$ copies of themselves, respectively.
\item $K_4''$ is $K_4$ with two adjacent edges doubled.
\item $K_1^\circ$ is, as in Section \ref{cycle}, a vertex with a loop.
\end{itemize}

\begin{thm}[Validity of the Circle Test] \label{T:good}
Let $\cG$ be the class of all groups, or all abelian groups, or any other subgroup-closed class that contains $\bbZ_3$.  For a finite graph $\G$, the following are equivalent:
\begin{enumerate}
\item[(i)] The circle test is valid for $\G$ with respect to $\cG$.
\item[(ii)] None of $C_3(3,3,2)$, $2C_4$, $K_4''$, or $W_4$ is a minor of $\G$.
\item[(iii)] Each block of $\G$ is isomorphic to a graph obtained by extrusion from one of $K_1^\circ$, $mK_2$ with $m \geq 1$, $C_3(m,2,2)$ with $m \geq 2$, or $K_4(m,m')$ with $m, m' \geq 1$.
\end{enumerate}
\end{thm}

The striking fact is that our forbidden minors depend entirely upon the presence of $\bbZ_3$ among the possible gain groups, even though there do exist bad graphs for all other odd cyclic groups.

The proof is by a series of auxiliary theorems that have some independent interest.  The first, Theorem \ref{T:fmcircle}, establishes the existence of the four forbidden minors and, going beyond that, shows that they are are good if $\bbZ_3$ is not a possible gain group.  Then none of the four is a forbidden minor, and indeed we know no forbidden minors.  The second auxiliary result, Theorem \ref{L:fmw4}, suggests a candidate forbidden minor if $\bbZ_{2k-1} \in \cG$, since $W_{2k}$ is bad if $\bbZ_{2k-1}$ is a possible gain group.  By the last auxiliary result, Theorem \ref{T:2sep}, a graph that has no $2C_4$ minor and is not obtained by an extrusion operation cannot have a 2-separation.  The proof of Theorem \ref{T:good} is then a short argument.

\subsection{General Methods}

To show a graph $\G$ is a forbidden minor we need to prove that $\G$ is bad and that every single-edge deletion or contraction of $\G$ is good.  To prove $\G$ is bad is fairly easy:  we produce a basis $B$ of balanced circles in an unbalanced gain graph based on $\G$.  
(By a slight abuse of terminology we say a basis \emph{$B$ is balanced} if every circle in $B$ is balanced.  We say $B$ \emph{implies balance} if every gain graph in which $B$ is balanced is itself balanced.)
It is much harder to prove a graph is good, because we have to treat every circle basis.

We can simplify this problem by the techniques of theta sums and theta summation.  The fundamental fact is this easy lemma:

\begin{lem}[Zaslavsky (1989)]  \label{L:theta}   
Let $C_1$ and $C_2$ be balanced circles in a gain graph.  If $C_1\cup C_2$ is a theta graph, then $C_1 + C_2$ is balanced.  
\end{lem}

Now suppose in a circle basis $B$ there are two circles, $C_1$ and $C_2$, whose union is a theta graph.  We call $C_1+C_2$ a \emph{theta sum}.  Then $B'$, obtained through replacing $C_1$ by $C_1+C_2$, is a circle basis and, moreover, if $C_1$ and $C_2$ are balanced, so is $C_1+C_2$.
Thus from one balanced circle basis we may obtain another that may be simpler.  We call this operation \emph{theta replacement}.  
A \emph{theta summation} is a sum of circles, $C_1+C_2+\cdots +C_k$, in which every $(C_1 + \cdots + C_i) + C_{i+1}$ is a theta sum.  Then if all $C_i$ are balanced, their sum is balanced.  

If $B$ is a basis for $Z_1(\G;\bbZ_2)$, then every edge belongs to at least one element of $B$.  If there is an edge $e$ that is in \emph{only} one $C \in B$, we call $e$ an \emph{improper edge} for $B$ and $B$ an \emph{improper basis} for $Z_1(\G;\bbZ_2)$.

\begin{lem} \label{L:proper}
Suppose $B$ is a binary circle basis of $\G$ with an improper edge $e$.  If $\G \setminus e$ is good, then $B$ implies balance in $\G$.
\end{lem}

\begin{proof}
$B \setminus e$ is a basis for $Z_1(\G \setminus e,\bbZ_2)$; hence $\G \setminus e$ is balanced.  Let $C'$ be a circle that contains $e$; we have to prove $C'$ is balanced.  Tutte's path theorem 
(Tutte (1956), (4.34)) 
tells us there is a chain of circles, $C=C_0, C_1, \ldots, C_l=C'$ such that each $C_{i-1} \cup C_i$ is a theta graph and $e \in C_i$ for $0<i<l$.  Each $C_{i-1} + C_i$ is in $\G \setminus e$ so is balanced.  Because $C$ is balanced it follows that $C'$ is balanced.  (This is essentially the proof that the balance-closure of a balanced set is balanced; 
see Zaslavsky (1989), Proposition 3.5.)
\end{proof}

A simple illustration of these methods is two proofs of a lemma about planar bases.

\begin{lem}  \label{L:faces}  
If $\G$ is a plane graph whose finite face boundaries are balanced circles in $\GG$, then $\GG$ is balanced.
\end{lem}

\begin{proof}[First proof]
Any circle is a theta summation of the boundaries of the faces it encloses.
\end{proof}

\begin{proof}[Second proof]
Any outer edge is improper with respect to the basis of face boundaries.  The lemma follows by induction.
\end{proof}

Theta sums lead to a valuable conclusion.  Suppose a balanced circle basis $B$ of $\GG$ contains a digon $e_1e_2$.  Using theta sums we can replace $e_2$ by $e_1$ in every other member of $B$; thus $B$ can be assumed to have the form $B' \cup \{e_1e_2\}$ where $B'$ is a balanced circle basis of $\G\setminus e_2$.  If $\G' = \G\setminus e_2$ is good with respect to a class $\cG$ (containing $G$), then $(\G',g',G)$ is balanced (where $g' = g\big|_{E\setminus e_2})$; it follows that $\GG$ is balanced.  
Let us say a circle basis $B$ satisfies the \emph{Digon Condition} with respect to $e_1e_2$ if it does not contain $e_1e_2$ and does not contain two circles whose sum is $e_1e_2$. 
We call a circle basis $B$ \emph{reducible} if, by one or more theta replacements, we can reduce the total length of all circles in $B$.

\begin{lem}[Digon Principle]  \label{L:digons}   
A graph $\G$ is good (with respect to a class $\cG$) if, for every gain graph $\GG$ with $G\in \cG$ and every proper circle basis $B$ of $\G$ that satisfies the Digon Condition with respect to every digon $e_1e_2$ such that $\G\setminus e_2$ is good, balance of $B$ implies balance of $\G$. 
$\hfill\blacksquare$
\end{lem}

\begin{lem} \label{L:mk2}
$mK_2$ is good for all $m \geq 1$ and any class $\cG$.
\end{lem}

\begin{proof}
Apply the Digon Principle.
\end{proof}

We need one more general lemma and one particular result.  A graph is \emph{inseparable} if any two edges lie in a common circle.  A \emph{block} of a graph is a maximal inseparable subgraph.

\begin{lem} \label{L:sep}
The Circle Test is valid for a graph if and only if it is valid for every block.
\end{lem}

\begin{proof}
Suppose there are gains on $\G$ such that some circle basis $B$ is balanced but $\GG$ is unbalanced.  Then there is an unbalanced block $\G'$.  The circles of $B$ that lie in $\G'$ form a circle basis for $\G'$ that is balanced.  Therefore, the Circle Test is invalid for $\G'$.

On the other hand, suppose there is a circle basis $B'$ for $\G'$ such that, for some gain map $g'$ on $\G'$, $B'$ is balanced but $\G'$ is not.  Extend $B'$ to a circle basis $B$ of $\G$ and extend $g'$ to $g$ on $\G$ by setting $g(e)=1$ for $e \notin E(\G')$.  Then $B$ is balanced but $\GG$ is not.
\end{proof}

This lemma lets us concentrate on inseparable graphs.  The next lemma establishes a fundamental example.

\begin{lem}  \label{L:k4mm}  
If $m, m' > 0$, then $K_4(m,m')$ is good for any class $\cG$.
\end{lem} 

\begin{proof} 
Let $e_{12}^i$ for $0<i\leq m$ and $e_{34}^j$ for $0<j\leq m'$ be the opposite sets of parallel edges.  (Of course, they are not multiple if $m$ or $m' = 1$.)
Our notation for triangles and quadrilaterals in $K_4(m,m')$ will be
\begin{gather*}
T^3_i = e_{12}^i e_{23} e_{31}, \qquad T^4_i = e_{12}^i e_{24} e_{41}, \\
T^1_j = e_{34}^j e_{41} e_{13}, \qquad T^2_j = e_{34}^j e_{42} e_{23}, \\
P_{ij} = e_{12}^i e_{23} e_{34}^j e_{41}, \qquad Q_{ij} = e_{12}^i e_{24} e_{43}^j e_{31}, \\
R = e_{13}e_{32}e_{24}e_{41}.
\end{gather*}
%

Deleting 
one of the parallel edges, say $e_{12}^m$, leaves a graph that can be assumed good by induction on $m$ if $m>1$ or is an extrusion of $(m'+2)K_2$, hence good, if $m = 1$.  Therefore the Digon Principle applies.  That is, we may assume we have a proper, irreducible circle basis $B$ that contains no two triangles $T^k_h$ for any fixed $k$, no two $P_{ij}$ and $P_{i'j'}$ or two $Q_{ij}$ and $Q_{i'j'}$ with $i=i'$ or $j=j'$, and no $P_{ij}$ or $Q_{ij}$ together with any $T^k_i$ ($k=3,4$) or $T^k_j$ ($k=1,2$).  It cannot include $R$ because $R$ is reducible in the presence of any triangle, and a basis cannot consist only of circles of even length.

We set up an auxiliary graph $A$ with vertex set $\{e_{12}^i, e_{34}^j : 0<i\leq m,\ 0<j\leq m'\}$ and edges those $P_{ij}$, $Q_{ij}$, and $T^k_h$ that belong to $B$.  $P_{ij}$ and $Q_{ij}$ have endpoints $e_{12}^i, e_{34}^j$.  $T^3_i$ and $T^4_i$ are half edges with endpoint $e_{12}^i$ and $T^1_j$ and $T^2_j$ are half edges with endpoint $e_{34}^j$.  (A half edge has only one endpoint; unlike a loop, it contributes one to the degree at that vertex.)  The edges of $A$ are thus quadrilaterals and triangles in $\G$.  The properties of $A$ and the reasons for them are:
\begin{enumerate}
\item Every vertex has degree at least 2.  (Reason:  $B$ is proper.)
\item No vertex is incident with two $P$ edges or two $Q$ edges or a $T$ edge together with a $P$ or $Q$ edge.  (Reason: irreducibility of $B$.)
\item No vertex is incident with more than two half edges.  (Each edge $e_{12}^i$ or $e_{34}^j$ lies in only two triangles.)
\item $A$ has $m+m'+1$ edges and $m+m'$ vertices.  (The number of edges is the dimension of the binary cycle space.)
\item $A$ contains no evenly even circle.  (Its edges would sum to zero in the binary cycle space.  ``Evenly even'' means the length is a multiple of 4.)
\item $A$ does not contain two oddly even circles, two vertices of type H, or one oddly even circle and one vertex of type H.  (An oddly even circle has even length that is not a multiple of 4.  A vertex has type H if it supports two half edges.  The reason: the sum of the edge quadrilaterals in an oddly even circle, or of the half-edge triangles at a vertex of type H, is $R$.)
\end{enumerate}

By (1)--(3) $A$ is regular of degree two, but by (4) it contains more edges than vertices.  That is possible only if $A$ has at least one vertex of type H.  If there is only one such vertex, then $A$ has a circle (because $m+m'>1$ so there are vertices besides the one of type H), hence it violates (5) or (6).  If there are two such vertices, $A$ violates (6).  
We conclude that no irreducible, proper basis $B$ exists; whence the gain graph is balanced by the Digon Principle and induction on $m$ and $m'$.
\end{proof}

\subsection{The four excluded minors} \label{4minors}

\begin{thm}  \label{T:fmcircle}  
Let $\cG$ be any subgroup-closed class containing $\bbZ_3$.  Then $\FM_0(\cG)$ contains $C_3(3,3,2)$, $2C_4$, $K_4''$, and $W_4$.  
However, if\/ $\bbZ_3 \notin \cG$, then all four graphs are good.
\end{thm}

To prove this we treat each graph in a separate proposition.

\begin{prop} \label{P:c3}  
The graph $C_3(3,3,2) \in \FM_0(\cG)$ if $\bbZ_3 \in \cG$, but it is good for $\cG$ if $\bbZ_3 \not\in \cG$.  
\end{prop}

\begin{proof}
By Lemma \ref{L:mk2}, $C_3(3,3,2)/e$ is good for every class $\cG$. 
Thus we need only consider deletion of one or more edges.

\begin{lem}  \label{L:c3(1)}  
If $m,m'\geq 1$, then $C_3(m,m',1)$ is good for every class $\cG$.  
\end{lem}

\begin{proof}  
By extrusion of $(m+m')K_2$.
\end{proof}

\begin{lem}  \label{L:c3(22)}  
If $m\geq 2$, then $C_3(m,2,2)$ is good for any class $\cG$.
\end{lem}

\begin{proof}  A basis that satisfies the digon requirements can have at most four members, but the cyclomatic number is $m+2$.  Thus the Digon Principle with induction implies goodness of $C_3(m,2,2)$ as long as $2C_3$ is good. 
$2C_3$ is the contraction of $K_4(2,1)$  by the edge $v_3v_4$, so the lemma follows from Lemma \ref{L:k4mm}.
\end{proof}

To complete the proof of the proposition we appeal to Example \ref{X:c3(332)} to show that $C_3(3,3,2)$ is bad if and only if $\bbZ_3 \in \cG$.  Moreover, deleting a digon edge we have $C_3(3,2,2)$ or $C_3(3,3,1)$, both of which are good. \end{proof}

\begin{prop}  \label{P:2c4}  
The graph $2C_4 \in \FM_0(\cG)$ if $\bbZ_3 \in \cG$, while it is good if $\bbZ_3 \not\in \cG$.  
\end{prop}

\begin{proof}  A deletion is an extrusion of $2C_3$, and a contraction is $2C_3$ with a loop; both are good by Lemma \ref{L:c3(22)}.  

Take a quadrilateral $Q$ in $2C_4$ and label its edges $e_1, e_2, e_3, e_4$ in cyclic order.  Let $f_i$ be parallel to $e_i$ and define $Q_i = Q + \{e_i,f_i\}$.  


We must demonstrate that $2C_4$ has unbalanced gains in $\bbZ_3$ that leave some circle basis $B_0$ balanced.  Let $B_0 = \{ Q,Q_1,Q_2,Q_3,Q_4 \}$.  
The gains are $g(e_i) =1$ and $g(f_i) = a$ in gain group 
$\langle a \mid a^3 = 1 \rangle$, the gains for all $f_i$ being calculated in a consistent direction around the quadrilateral $f_1f_2f_3f_4$.  

Conversely, suppose unbalanced gains in a group $G$ leave $B_0$ balanced.  We may switch so all $g(e_i) = 1$.  It is easy to deduce that the gains are those just described:  that is, if $G\not\geq \bbZ_3$ and $B_0$ is balanced, then $(2C_4,g,G)$ is balanced.

Finally, suppose $\bbZ_3 \not\in \cG$.  By Lemma \ref{L:digons} and the goodness of $2C_4 \setminus f_4$ we need only consider bases $B$ that satisfy the Digon Condition.  We encode quadrilaterals as binary sequences, 0 or 1 in position $i$ corresponding to edge $e_i$ or $f_i$, and we assume $0000\in B$.  For $B$ to be a cycle basis it must contain a sequence of odd weight; thus, say $1110 \in B$.  If a second even sequence and a second odd sequence belong to $B$, $B$ cannot contain the necessary five elements.  
Therefore either $B$ has no sequence of even weight other than 0000, in which case $B=B_0$, or $B$ contains only the one odd sequence $1110$, in which case $B = \{0000,1110,1001,0101,0011\}$.  Assuming by adequate switching that $g(e_i) = 1$, we deduce that $g(f_i) = a$ for $i = 1,2,3$ and $g(f_4) = a^{-1}$ where $a^3 = 1$.  Again the gain group after switching is $\langle a\mid a^3 = 1\rangle$.  

Thus in every case if an unbalanced $2C_4$ has a balanced circle basis, the gain group contains $\bbZ_3$.  
\end{proof}

\begin{prop} \label{P:k4dd}
The graph $K_4'' \in \FM_0(\cG)$ if $\bbZ_3 \in \cG$, but it is good if $\bbZ_3 \notin \cG$.
\end{prop}

\begin{proof}
Deleting an edge gives a good graph for any class $\cG$, by Lemmas \ref{L:c3(1)}, \ref{L:c3(22)}, and \ref{L:k4mm} and extrusion.  Contraction also gives a good graph.

\begin{figure}
\vskip 2.5in{}
\includegraphics{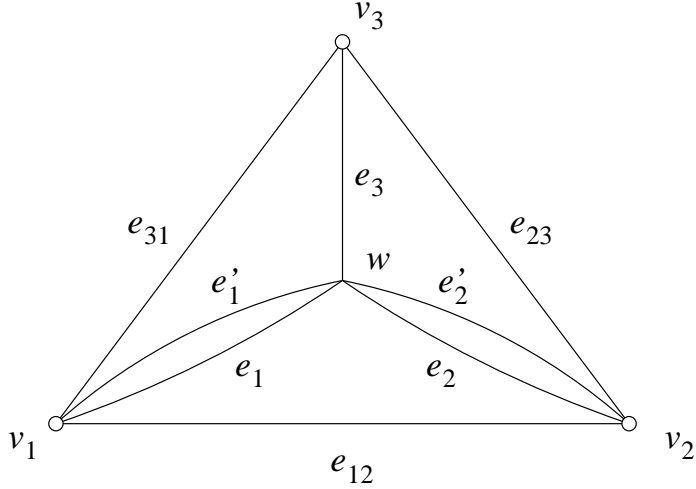}
\caption{The edge labels of $K_4''$.} \label{F:k4dd}
\end{figure}

For $K_4''$ itself (see Figure \ref{F:k4dd}) we take the binary cycle basis $B = \{Q_1,Q_2,T_{12},Q_3^L,Q_3^R \}$ where
\begin{gather*}
Q_i = e_{i3} e_{12} e_{3-i}' e_3, \qquad 
T_{12} = e_{12} e_1 e_2, \\
Q_3^L = e_1' e_2 e_{23} e_{31}, \qquad Q_3^R = e_1 e_2' e_{23} e_{31}.
\end{gather*}
Assuming each of these has identity gain, we compute the gain group.  By switching we may assume $g(e_i) = 1$ for $i=1,2,3$.  Orienting $e_i' $ from $w$ to $v_i$, then
\begin{equation*}
\begin{aligned}
T_{12} &\implies g(e_{12}) = 1, \\
Q_1 &\implies g(e_{31})g(e_{12})g(e_2')\inv = 1 \implies g(e_{31}) = g(e_2') = a, \\
Q_2 &\implies g(e_{32}) = g(e_1') = b, \\
Q_3^L &\implies a = g(e_{31}) = g(e_{32})g(e_1') = b^2, \\
Q_3^R &\implies b = g(e_{32}) = g(e_{31})g(e_2') = a^2,
\end{aligned}
\end{equation*}
from which it follows that $a^3 = 1$ and $b=a^2$.  These are the only relations; thus $K_4''$ is bad if $\bbZ_3 \in \cG$ and good otherwise.
\end{proof}

\begin{prop}  \label{P:w4}  
The wheel $W_4 \in \FM_0(\cG)$ if $\cG$ contains $\bbZ_3$ but it is good otherwise.
\end{prop}

The first lemma is the case of $W_4$ in a valuable result that shows there is a bad graph for every odd cyclic gain group.

\begin{thm} \label{L:fmw4}  
For each $k \geq 2$, $W_{2k}$ and $2C_{2k}$ are bad if any of $\bbZ_3, \bbZ_5, \ldots, \bbZ_{2k-1} \in \cG$.  
\end{thm}

\begin{proof}  
The $2k$ Hamiltonian circles of $W_{2k}$ constitute a basis $B_{2k+1}$ for $Z_1(W_{2k},\bbZ_2)$.  
Suppose they are all balanced in a gain graph $(W_{2k},g,G)$.  
Write $H_i = wv_{i}v_{i+1}\cdots v_{i-1}w$ (subscripts taken modulo $2k$).  By switching we may assume all $g(wv_i) = 1$.  Letting $g_i = g(v_{i-1}v_i)$, we have 
\begin{equation*}
g_{i+1}g_{i+2}\cdots g_{i-1} = g(H_i) = 1 \text{ for }i = 1,\ldots,2k.
\end{equation*}
It follows that $g_1 = g_2 = \ldots = g_{2k}$ and $g_1^{2k-1} = 1$.  Thus the gain group may be $\bbZ_{2k-1}$ with $g_1$ as a generator. Clearly then all $H_i$ are balanced while the gain graph is unbalanced.  It follows that $W_{2k}$ is bad if $\bbZ_{2k-1} \in \cG$.  
Since $W_{n-1}$ is a minor of $W_n$, the whole result for wheels follows.

For $2C_{2k}$ we generalize the relevant part of the proof of Proposition \ref{P:2c4}.  Call the edges $e_i, f_i$ for $i=1,2,\ldots,2k$, $e_i$ and $f_i$ being parallel, and let $C=e_1e_2\cdots e_{2k}$, $D=f_1 f_2 \cdots f_{2k}$, and $D_i = D + \{e_i,f_i\}$ (this is set sum).  As balanced circle basis take $B = \{C, D_1, \ldots, D_{2k} \}$.  By switching assume all $g(e_i)=1$.  Then all the gains $g(f_i)=a$ with $a^{2k-1} = 1$.  Obviously, $\bbZ_{2k-1}$ is a possible gain group.
\end{proof}

\begin{proof}[Proof of Proposition \ref{P:w4}]  
We proved $W_4$ is bad if $\bbZ_3 \in \cG$.  Contracting an edge in $W_4$ gives the graph $K_4(2,1)$ of Lemma \ref{L:k4mm} or $K_4''$ with a particular simple edge deleted, 
both good (the latter by extrusion from $5K_2$ or $C_3(2,2,2)$).  
Deleting an edge gives a subdivision of $K_4$, which is good by Lemma \ref{L:k4mm}, or of $C_3(2,2,1)$; the latter is good by Lemma \ref{L:c3(1)}.


Now we assume that $\bbZ_3 \notin \cG$ and show that $W_4$ is good.  In fact, the Hamiltonian basis $B_5$ is the only one that can fail to imply balance of $(W_4,g,G)$, and we showed in the proof of Theorem \ref{L:fmw4} that it implies balance if $G\ngeq\bbZ_3$.

\begin{lem} \label{L:w4bal}
Let $B$ be a circle basis of $W_4$ other than $B_5$.  If $(W_4,g,G)$ is a gain graph in which every circle in $B$ is balanced, then it is balanced.
\end{lem}

\begin{proof}
The circle basis $B_3$, the set of triangles, implies balance for all gain groups by Lemma \ref{L:faces}.   We shall prove that every circle basis except for $B_5$ and of course $B_3$ is reducible to $B_3$.
Since balance of $B$ implies balance of any basis obtained by theta replacement, $B$ is balanced$\implies$$B_3$ is balanced$\implies$the gain graph is balanced.

Circles whose appearance together in a basis implies reducibility are:

\begin{enumerate}
\item[(a)] A Hamiltonian circle $H_i$ and a triangle $T_j = wv_{j-1}v_{j}w$ with $j = i-1, i, i+1$.  ($H_i$ was defined at Theorem \ref{L:fmw4}.  Subscripts are modulo 4.)
\item[(b)] $H_i$ and a quadrilateral $Q_j = wv_{j-1}v_jv_{j+1}w$ for $j=i-1$ and $i$.   (We call $Q_{i-1}$ and $Q_{i}$ \emph{consecutive}.)
\item[(c)] $T_i$ and $Q_{i-1}$ or $Q_{i}$.
\end{enumerate}

Consider an irreducible basis that contains $H_i$ and $T_{i-2}$.  The remaining basis elements must be $Q_{i-1}$ and $Q_{i}$, by (a) and (b).  But $H_i + T_{i-2} + Q_{i-1} + Q_{i} = 0$, so there is no such basis.

Consider an irreducible basis that contains a Hamiltonian circle but is not $B_5$.  No triangle can be in the basis, nor can the rim quadrilateral.  Since each Hamiltonian circle allows only two consecutive nonrim quadrilaterals, one can easily see that two Hamiltonian circles allow at most one quadrilateral, and three allow none.  Because we cannot have the required four circles, no such irreducible basis exists.

Consider finally an irreducible basis $B$ that contains a triangle but is not $B_3$.  Suppose first that the rim, $R$, is in $B$.  If the rest of $B$ is three triangles, $R$ can be replaced by the fourth triangle, reducing $B$ to $B_3$.  Otherwise, some $Q_i$ is in $B$; by theta replacement we can substitute $Q_{i-2} = Q_i + R$ for $R$.  Thus we may assume $R$ is not a basis element: $B$ consists of triangles and at least one $Q_i$.  (It cannot consist only of quadrilaterals because they will not generate odd circles.)  Now reasoning as with a basis consisting of Hamiltonian circles and quadrilaterals, we conclude that no irreducible basis exists.
\end{proof}

Clearly, Lemma \ref{L:w4bal} implies the proposition.
\end{proof}

We suggest that $W_{2k}$ and $2C_{2k} \in \FM_0(\cG)$ if (and only if) $\cG$ contains $\bbZ_m$ for $m$ a factor of $2k-1$.

\subsection{2-Separations}\label{2sum}

The \emph{parallel connection} or \emph{edge amalgamation} of two graphs is obtained by assuming the graphs are disjoint and identifying an edge, say $e_1$, in the first with an edge, say $e_2$, in the second.  The 2-\emph{sum} of the graphs is the parallel connection with the identified edge deleted.  (One may require the identified edges to be links.  That will not affect our discussion.)  
A natural question is whether the validity of any of the tests, given a class of possible gain groups, is preserved under parallel connection or 2-summation.  If that were true, every forbidden minor would be 3-connected or very small.
The example of $2C_4$ shows that the validity of the Circle Test is not preserved under 2-summation, and therefore not under parallel connection.  However, $2C_4$ is perhaps too special; here is a possibly more representative modification in which the 2-separation is unique.

\begin{exam}\label{X:2sum} {\rm 
Take two disjoint copies of $K_4$, $\G_1$ and $\G_2$.  In $\G_i$ take two nonadjacent edges, called $e_i$ and $f_i$.  Take the 2-sum along $e_1$ and $e_2$ and contract $f_1$ and $f_2$.  Now, $K_4$ is good for the Circle Test for any gain group (Lemma \ref{L:k4mm}) but the 2-sum is bad because the contraction is $2C_4$, which is bad for the gain group $\bbZ_3$ (Proposition \ref{P:2c4}).
}\end{exam}

On the other hand, consider these 2-separable graphs:

\begin{exam}\label{X:tripart} {\rm 
$K_{11p}(m;m_1,\ldots,m_p)$ is obtained by taking the complete tripartite graph $K_{11p}$ with vertices $v$, $w$, and $x_1,\ldots,x_p$, replacing one of the two edges at each vertex $x_i$ by $m_i$ copies of itself, and replacing $vw$ by $m$ copies of itself.  For any $p, m, m_1,\ldots, m_p > 0$, this graph is good because it is obtained by extrusion from $(m+m_1+\cdots+m_p)K_2$.
}
\end{exam}

\begin{exam}\label{X:tripartd} {\rm 
Let $m'_1, m'_2 \geq 2$.  $K_{1,1,p+1}(m;m_1,\ldots,m_p;m'_1,m'_2)$ is obtained by taking $K_{11p}(m;m_1,\ldots,m_p)$ and adding a $vw$-path of length at least two, two of whose edges are replaced by $m'_1$ and $m'_2$ copies of themselves.  
This graph is good if $m'_1=m'_2=2$ and bad otherwise (unless $m=p=m_1=1$) because it is obtained by extrusion from $C_3(m+m_1+\cdots+m_p,m'_1,m'_2)$.
}
\end{exam}

\begin{exam}\label{X:tripartk} {\rm 
$K'_{11p}(m;m_1,\ldots,m_p)$ is obtained by subdividing one of the $vw$ edges in
$K_{11p}(m;m_1,\ldots,m_p)$ into a three-edge path $P_3$ and taking the parallel connection with $K_4$ along the middle edge of the path.
It is good because it is obtained by extrusion from $K_4(m+m_1+\cdots+m_p)$.
}
\end{exam}

All these examples, good or bad, are obtained by extrusion.  With that observation as guide we prove that $2C_4$ is the only 2-separable forbidden minor for the Circle Test with respect to a class $\cG$ that includes $\bbZ_3$.
The major part of the proof is a theorem that is not directly connected with the Circle Test.  

\begin{thm}  \label{T:2sep}  
Let $\Gamma$ be a finite, inseparable, extrusion-irreducible graph of which $2C_4$ is not a minor.  Then $\Gamma$ has no 2-separation.
\end{thm}

\begin{proof}
This requires some definitions and lemmas.  
Let $\Gamma$ be an arbitrary inseparable graph and $u,v\in V(\Gamma)$.  A \emph{bridge of} $\{u,v\}$ is a maximal subgraph $\Delta$ of $\Gamma$ such that any two elements of $\Delta$, whether vertices or edges, lie in a common walk that is internally disjoint from $\{u,v\}$.  For instance, a $\{u,v\}$-bridge may have just a single edge.  We classify bridges into three sorts.  Let $2P_{2uv}$ denote a graph that consists of $2P_2$ ($P_2$ being a path of length two) with $u$ and $v$ as its endpoints.  A non-edge bridge of $\{u,v\}$ has \emph{type II} if it has $2P_{2uv}$ as a minor and \emph{type I} otherwise.  A \emph{2-bridge} of $\Gamma$ is any bridge of any pair of vertices; each 2-bridge therefore has type I or type II or is an edge.  

We write $\Gamma\geq M$ to mean that $\Gamma$ has a minor isomorphic to $M$.

\begin{lem}  \label{L:2c4bridge}  
If $\Gamma\not\geq 2C_4$, then at most one bridge of any pair $\{u,v\}$ has type II.
\end{lem}

\begin{proof}  Obvious.
\end{proof}

A \emph{separating vertex} of a bridge $\Delta$ of $\{u,v\}$ is a vertex $w$ such that $\Delta\setminus w$ is disconnected.  It is clear that $w$ cannot be $u$ or $v$.

\begin{lem}  \label{L:sepbridge}  
Suppose $\Gamma$ is inseparable.  Any 2-bridge of type I has a separating vertex.
\end{lem}

\begin{proof}  
Suppose a bridge $\Delta$ of $\{u,v\}$ has no separating vertex.  By Menger's theorem there exist internally disjoint $uv$-paths $P$ and $Q$.  Since $P$ and $Q$ are in the same bridge, there is a path $R$ in $\Delta\setminus \{u,v\}$, joining $P$ to $Q$ and internally disjoint from both.  Then $P\cup Q\cup R$ contracts to $2P_{2uv}$.  
\end{proof}

\begin{lem}  \label{L:type1}  
If $\Gamma$ is finite, inseparable, and extrusion-minimal, then it has no 2-bridge of type I.
\end{lem}

\begin{proof}  Take a type I bridge $\Delta$ of $\{u,v\}$; thus $\Delta$ has a separating vertex $w$ which splits it into $\Delta_u \owns u$ and $\Delta_v \owns v$.  That $\Delta$ has type I means that one of $\Delta_u$ and $\Delta_v$ is a path (hence an edge) while in the other $w$ has degree at least two.  Say $\Delta_u$ is the path.  Amongst the bridges of $\{w,v\}$ is one that contains all the
bridges of $\{u,v\}$ except $\Delta$.  All other bridges $\Delta'$ of $\{w,v\}$ are contained in $\Delta$ and are edge bridges or have type I.  If $\Delta'$ is not an edge, it is a 2-bridge of type I that is properly contained in $\Delta$.  

Now suppose $\Delta$ in the preceding discussion to be a \emph{minimal} 2-bridge of type I.  Then any $\Delta'$ is an edge.  Therefore $\Delta$ is $P_2$ with vertices $u,w,v$ and with edge $wv$ replaced by a multiple edge.  Thus, $\Gamma$ is obtained by extruding $w$ from $v$ in the contraction $\Gamma/uw$, contrary to hypothesis.  It follows that no 2-bridge of type I can exist.  
\end{proof}

To complete the proof of Theorem \ref{T:2sep}, take any vertex pair $\{u,v\}$ in $\Gamma$.  Every bridge but one is an edge.  Therefore, $\{u,v\}$ cannot separate $\Gamma$.
\end{proof}


\subsection{The End}

\begin{proof} [Proof of Theorem \ref{T:good}]  
We proved (i)$\implies$(ii) in Theorem \ref{T:fmcircle}.  We know (iii)$\implies$(i) by Lemma \ref{L:extrude} and the lemmas that imply all the graphs listed in (iii) are good.

As for (ii)$\implies$(iii), it suffices to prove it for graphs that are inseparable (since all the forbidden minors in (ii) are inseparable) and extrusion-irreducible.  To see the latter, suppose an inseparable graph $\Gamma$ is obtained by repeated extrusion from an extrusion-irreducible graph $\Gamma_0$, which of course satisfies (ii) and is inseparable.  Then $\Gamma_0$ is one of the list in (iii) so $\Gamma$ is obtained as in (iii).

Thus, let $\Gamma$ be an inseparable, extrusion-irreducible graph that satisfies (ii).  By Theorem \ref{T:2sep}, $\Gamma$ is 3-connected or has order at most three.  Any 3-connected graph of order 5 or more has $W_4$ as a minor.  A 3-connected graph of order 4 is $K_4$ with, possibly, multiple edges; this must be $K_4(m,m')$ because $K''_4$ is excluded.  An inseparable graph of order 3 is
$C_3(m_1,m_2,m_3)$ with, say, $m_1\geq m_2 \geq m_3 > 0$.  By extrusion irreducibility, $m_3 \geq 2$.
By exclusion of $C_3(3,3,2)$, $m_2 = m_3 = 2$; thus the graph is $C_3(m,2,2)$.  Graphs of order 1 or 2 are good.  Thus the theorem is proved.  
\end{proof}


\end{document}